\documentclass[12pt,twoside]{article}
\usepackage[english]{babel}
\usepackage[latin1]{inputenc}
\usepackage{amsmath}
\usepackage{amssymb,amsfonts}
\usepackage{graphicx}
\usepackage{times,amssymb,amscd}

\newcommand{\bG}{\mathbf{G}}
\newcommand{\bH}{\mathbf{H}}

\newcommand{\bL}{\mathbf{L}}

\newcommand{\bR}{\mathbf{R}}

\newcommand{\bS}{\mathbf{S}}

\newcommand{\bg}{\mathbf{g}}

\newcommand{\bI}{\mathbf{I}}

\newcommand{\cD}{\mathcal{D}}

\newcommand{\cP}{\mathcal{P}}

\newcommand{\cB}{\mathcal{B}}

\newcommand{\cA}{\mathcal{A}}

\newcommand{\SPH}{\bS^3}
\newcommand{\HYP}{\bH^3}
\newcommand{\SXR}{\bS^2\!\times\!\bR}
\newcommand{\HXR}{\bH^2\!\times\!\bR}
\newcommand{\SLR}{\widetilde{\bS\bL_2\bR}}
\newcommand{\NIL}{\mathbf{Nil}}
\newcommand{\SOL}{\mathbf{Sol}}

\begin{document}
\pagestyle{myheadings}
\markboth{\centerline{B.~Schultz,~J. Szirmai }}
{Densest geodesic ball packings...}
\title
{Densest geodesic ball packings to $\mathbf{S}^2\!\times\!\mathbf{R}$ space groups generated by screw motions \footnote{AMS Classification 2000: 52C17, 52C22, 53A35, 51M20}}

\author{ B.~Schultz, ~ J. Szirmai\footnote{
 E-mail:
schultzb@math.bme.hu,
szirmai@math.bme.hu
}\\
\normalsize Budapest University of Technology and \\
\normalsize Economics Institute of Mathematics, \\
\normalsize Department of Geometry\\
}


\maketitle
\begin{abstract}

In this paper we study the locally optimal geodesic ball packings with equal balls to the $\SXR$ space groups having rotation point groups
and their generators are screw motions.
We determine and visualize the densest simply transitive geodesic ball arrangements for the above space groups, moreover we compute their optimal
densities and radii. The densest packing is derived from the $\SXR$ space group $\mathbf{3qe.~I.~3}$ with packing density $\approx 0.7278$.

E. {Moln\'ar} has shown in \cite{M97}, that the Thurston geometries
have an unified interpretation in the real projective 3-sphere $\mathcal{PS}^3$.
In our work we shall use this projective model of $\SXR$ geometry.
\end{abstract}

\newtheorem{Theorem}{Theorem}[section]
\newtheorem{corollary}[Theorem]{Corollary}
\newtheorem{lemma}[Theorem]{Lemma}
\newtheorem{exmple}[Theorem]{Example}
\newtheorem{definition}[Theorem]{Definition}
\newtheorem{Remark}[Theorem]{Remark}
\newtheorem{proposition}[Theorem]{Proposition}
\newenvironment{remark}{\begin{rmrk}\normalfont}{\end{rmrk}}
\newenvironment{example}{\begin{exmple}\normalfont}{\end{exmple}}
\newenvironment{acknowledgement}{Acknowledgement}



\section{Introduction}\label{sec_1}

The packing of spheres has been studied by mathematicians for centuries. Finding the densest packing
of balls in the Euclidean space is called the Kepler-problem or conjecture, which was solved
by Thomas Hales, who proved, with the use of computers, that Keplers original guess was correct \cite{H}.

In the $3$-dimensional hyperbolic $\HYP$ and spherical space $\SPH$ the problem of the ball packing was widely investigated
in a lot of papers (e.g. \cite{BF64}), but there are several open questions in this topic for example related to the
horoball and hyperball packings in $\HYP$ (see \cite{KSz}, \cite{Sz11-3}, \cite{Sz12}).

In \cite{Sz14-1} the second author has extended the problem of finding the densest geodesic and translation
ball (or sphere) packing for the other $3$-dimensional homogeneous geometries (Thurston geometries)
$\SXR,~\HXR,~$ $\SLR,~\NIL,~\SOL$,  \cite{Sz07-2}, \cite{Sz10}, \cite{Sz11-2}, \cite{Sz11-4}, \cite{Sz11}.

In this paper we consider the $\SXR$ geometry which can be derived by the direct product of the spherical plane $\bS^2$ and the real line $\bR$.
In \cite{F01} J.~Z.~Farkas has classified and given the complete list of the space groups in $\SXR$.
The $\SXR$ manifolds up to similarity and diffeomorphism were classified by E.~Moln\'ar and J.~Z.~Farkas in \cite{FM01}.
In \cite{Sz11-1} the geodesic balls and their volumes were studied, moreover, we have introduced
the notion of the geodesic ball packing
and its density and determined the densest simply and multiply
transitive geodesic ball packings for generalized Coxeter space groups of $\SXR$, respectively.
Among these groups the density of the densest packing is $\approx 0.8245$.

In \cite{Sz11-2} we have studied the simply transitive locally optimal ball packings to the $\SXR$ space groups having Coxeter point groups
and at least one of the generators is a non-trivial glide reflection.
We have determined the densest simply transitive geodesic ball arrangements for the above space groups. Moreover, their optimal
densities and radii were determined.
The density of the densest packing of the above space groups is $\approx 0.8041$.

Moreover, a candidate of the densest geodesic ball
packing is described in \cite{Sz14-1}. In the Thurston geometries the greatest known density was $\approx 0.8533$
that is not realized by a packing with {\it equal balls} of the hyperbolic
space $\HYP$. However, that is attained, e.g., by a {\it horoball packing} of
$\overline{\bH}^3$ where the ideal centres of horoballs lie on the
absolute figure of $\overline{\bH}^3$ inducing the regular ideal
simplex tiling $(3,3,6)$ by its Coxeter-Schl\"afli symbol.
In \cite{Sz14-1} we have presented a geodesic ball packing in the $\SXR$ geometry
whose density is $\approx 0.8776$.

In this paper we study the locally optimal ball packings to the $\SXR$ space groups belonging to the {\it screw motion groups}, i.e.,
the generators \linebreak$\bg_i \ (i=1,2,\dots m)$ of its point group $\Gamma_0$
are rotations, and either any translation parts are zero or at least one of the possible translation parts of
the above generators differs from zero (see \cite{Sz11-2}).
We determine and visualize the densest simply transitive geodesic ball arrangements for the above space groups, moreover we compute their optimal
densities and radii.
\section{The structure of the $\SXR$ geometry}
Now, we shall discuss the simply transitive ball packings to a given space group. But let us start first with the necessary concepts.

$\SXR$ geometry can be derived by the direct product of the spherical plane $\bS^2$ and the real line $\bR$.
The points in $\SXR$ geometry are described by $(P,p)$ where $P\in \bS^2$ and $p\in \bR$.

The isometry group $Isom(\SXR)$ of $\SXR$ can be derived by the direct product
of the isometry group of the spherical plane $Isom(\bS^2)$ and the isometry group of the real line $Isom(\bR)$.
The structure of an isometry group $\Gamma \subset Isom(\SXR)$  is the following: $\Gamma:=\{(A_1 \times \rho_1), \dots (A_n \times \rho_n) \}$, where
$A_i \times \rho_i:=A_i \times (R_i,r_i):=(g_i,r_i)$, $i \in \{ 1,2, \dots n \}$, and $A_i \in Isom(\bS^2)$, $R_i$ is either the identity map
$\mathbf{1_R}$ of $\bR$ or the point reflection $\overline{\mathbf{1}}_{\mathbf{R}}$. $g_i:=A_i \times R_i$ is called the linear part of the transformation
$(A_i \times \rho_i)$ and $r_i$ is its translation part.

The multiplication formula is the following:
\begin{equation}
(A_1 \times R_1,r_1) \circ (A_2 \times R_2,r_2)=((A_1A_2 \times R_1R_2,r_1R_2+r_2). \tag{2.1}
\end{equation}
A group of isometries $\Gamma \subset Isom(\SXR)$ is called {\it space group} if the linear parts form a finite group $\Gamma_0$
called the point group of
$\Gamma$. Moreover, the translation parts to the identity of this point group are required to form a one dimensional lattice $L_{\Gamma}$ of $\bR$.
It can be proved that the space group $\Gamma$ exactly described above has a compact fundamental domain $\mathcal{F}_\Gamma$.
We characterize the spherical plane groups by the {\it Macbeath-signature} (see \cite{M}, \cite{T}).

{\it In this paper we deal with a class of the $\SXR$ space groups {\bf 1q.~I.~1}; 
{\bf 1q.~I.~2}; {\bf 3q.~I.~1}; {\bf 3q.~I.~2}; {\bf 3qe.~I.~3}; {\bf 8.~I.~1}; {\bf 8.~I.~2}; {\bf 9.~I.~1}; {\bf 9.~I.~2}; {\bf 10.~I.~1};
(with a natural parameter $q \ge 2$, see \cite{F01}) where each of them belongs to the {\it screw motion groups}}.
\subsection{Geodesic curves and balls in $\SXR$ space}
In \cite{Sz11-1} and \cite{Sz11-2} we have described the equation system of the geodesic curve and so the geodesic sphere:
\begin{equation}
  \begin{gathered}
   x(\tau)=e^{\tau \sin{v}} \cos{(\tau \cos{v})}, \
   y(\tau)=e^{\tau \sin{v}} \sin{(\tau \cos{v})} \cos{u}, \\
   z(\tau)=e^{\tau \sin{v}} \sin{(\tau \cos{v})} \sin{u},\
   -\pi < u \le \pi,\ \ -\frac{\pi}{2}\le v \le \frac{\pi}{2} \tag{2.2}
  \end{gathered}
\end{equation}
with radius $\rho=\tau \ge 0$ and centre $x(0)=1,~y(0)=0,~z(0)=0,$ and longitude $u$, altitude $v$, as geographical coordinates.

In \cite{Sz11-1} we have proved that a geodesic sphere $S(\rho)$ in $\SXR$ space is a simply connected surface
in $\mathbf{E}^3$ if and only if $\rho \in [0,\pi)$. Namely, if $\rho\ge\pi$ then there is at least one $v \in
[-\frac{\pi}{2},\frac{\pi}{2}]$ so that $y(\tau,v)=z(\tau,v)=0$, i.e., selfintersection would occur.
Thus we obtain the following
\begin{proposition}
The geodesic sphere and ball with radius $\rho$ exists in the $\SXR$ space if and only if $\rho \in [0,\pi).$
\end{proposition}
We have obtained (see \cite{Sz11-1}) the volume formula of the geodesic ball $B(\rho)$ of radius $\rho$ by
the metric tensor $g_{ij}$, by the Jacobian of (2.2)
and a careful numerical Maple computation for given $\rho$ by the following integral:
\begin{Theorem}
\begin{equation}
\begin{gathered}
Vol(B(\rho))=\int_{V} \frac{1}{(x^2+y^2+z^2)^{3/2}}\mathrm{d}x ~ \mathrm{d}y ~ \mathrm{d}z = \\ = \int_{0}^{\rho} \int_{-\frac{\pi}{2}}^{\frac{\pi}{2}}
\int_{-\pi}^{\pi}
|\tau \cdot \sin(\cos(v)\tau)| ~ \mathrm{d} u \ \mathrm{d}v \ \mathrm{d}\tau =
2 \pi \int_{0}^{\rho} \int_{-\frac{\pi}{2}}^{\frac{\pi}{2}} |\tau \cdot \sin(\cos(v)\tau)| ~ \mathrm{d} v \ \mathrm{d}\tau. \tag{2.3}
\end{gathered}
\end{equation}
\end{Theorem}
\subsection{On fundamental domains}
A type of the fundamental domain of a studied space group can be combined as a fundamental domain of the corresponding spherical
group with a part of a real
line segment. This domain is called $\SXR$ {\it prism} (see \cite{Sz11-1}).
{\it This notion will be important to compute the volume of the Dirichlet-Voronoi cell of a given space group because their volumes are equal and
the volume of a $\SXR$ prism can be calculated by Theorem 2.3.}
\begin{figure}[ht]
\centering
\includegraphics[width=6.5cm]{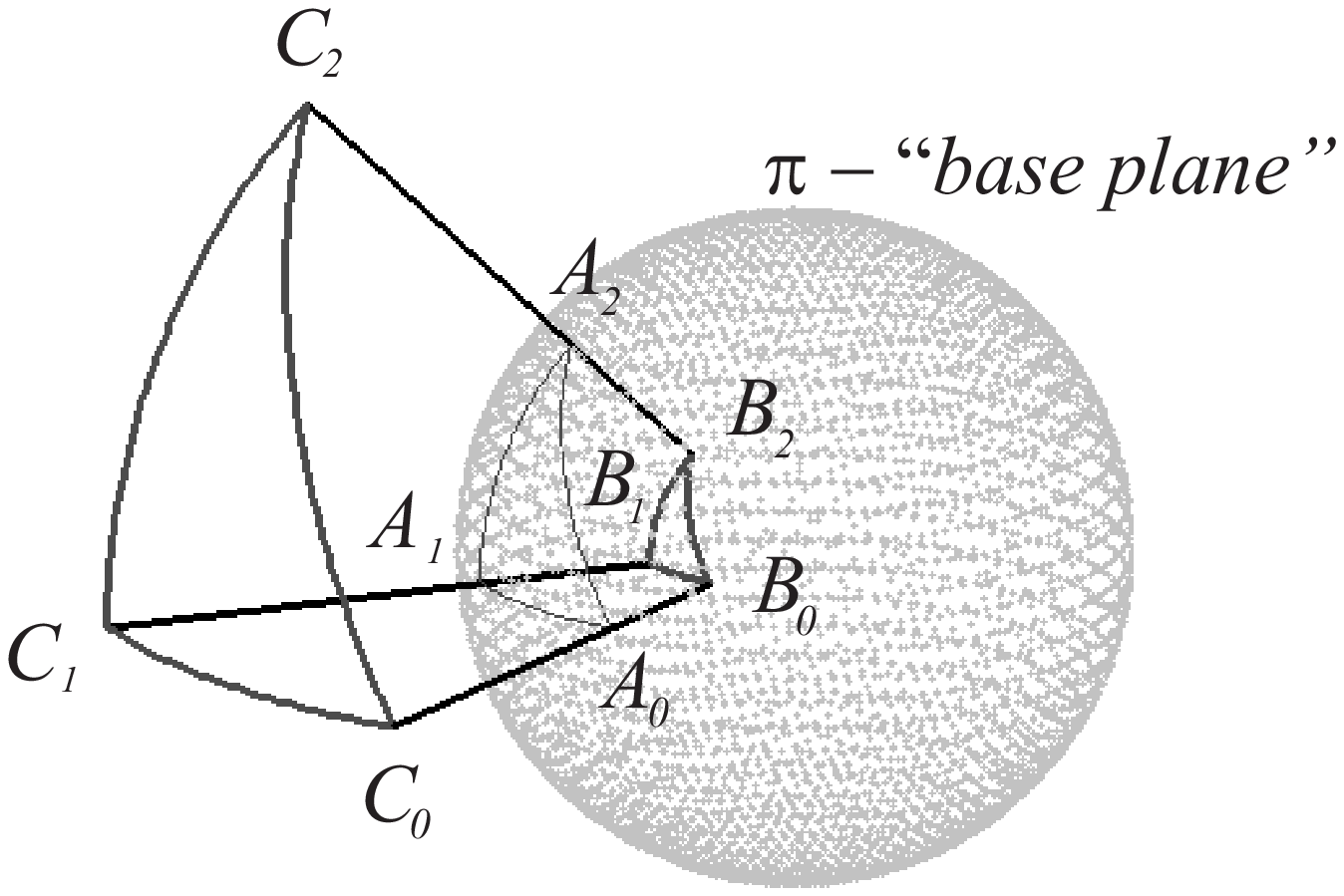} \includegraphics[width=5cm]{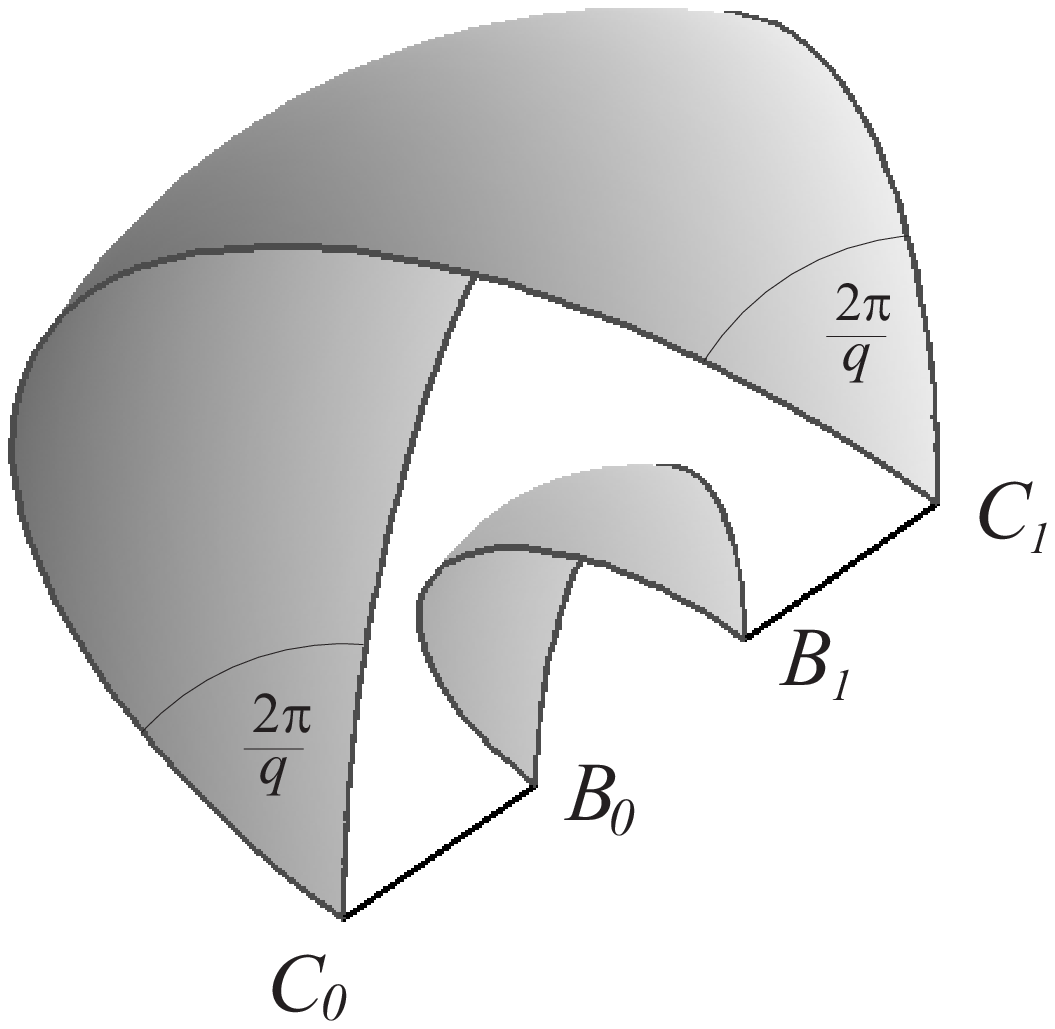}
\caption{Prism-like fundamental domains}
\label{pic:prismfund}
\end{figure}
The $p$-gonal faces of a prism called cover-faces, and the other faces
are the side-faces. The midpoints of the side edges
form a "spherical plane" denoted by $\Pi$. It can be assumed that the plane $\Pi$ is the {\it base plane}
in our coordinate system (see (2.2)) i.e. the fibre coordinate $t=0$.

From \cite{Sz11-1} we recall
\begin{Theorem}
The volume of a $\SXR$ trigonal prism $\mathcal{P}_{B_0B_1B_2C_0C_1C_2}$ and of a digonal prism
$\mathcal{P}_{B_0B_1C_0C_1}$ in $\SXR$ (see Fig.~1.a-b) can be computed by the following formula:
\begin{equation}
Vol(\mathcal{P})=Vol(\mathcal{A}) \cdot h \tag{2.4}
\end{equation}
where $Vol(\mathcal{A})$ is the area of the spherical triangle $A_0A_1A_2$ or digon $A_0A_1$ in the base plane $\Pi$ with fibre
coordinate $t=0$, and
$h=B_0C_0$ is the height of the prism.
\label{ttt_prismvol}
\end{Theorem}
\section{Geodesic ball packings \\ under discrete isometry groups}
A $\SXR$ space group $\Gamma$ has a compact fundamental domain.
Usually the shape of the fundamental domain of a group of $\bS^2$ is not
determined uniquely but the area of the domain is finite and unique
by its combinatorial measure. Thus the shape of the fundamental domain of a crystallographic group of $\SXR$ is not unique, as well.

In the following let $\Gamma$ be a fixed by {\it screw motions generated} space group of $\SXR$. We
will denote by $d(X,Y)$ the distance of two points $X$, $Y$ by definition (2.2).
\begin{definition}
We say that the point set
$$
\cD(K)=\{X\in\SXR\,:\,d(K,X)\leq d(K^\bg,X)\text{ for all }\bg\in\ \Gamma\}
$$
is the \emph{Dirichlet--Voronoi cell} (D-V~cell) to $\Gamma$ around the kernel
point $K\in\SXR$.
\end{definition}
\begin{definition}
We say that
$$
\Gamma_X=\{\bg\in\Gamma\,:\,X^\bg=X\}
$$
is the \emph{stabilizer subgroup} of $X\in\SXR$ in $\Gamma$.
\end{definition}
\begin{definition}
Assume that the stabilizer $\Gamma_K=\bI$ i.e. $\Gamma$ acts simply transitively on
the orbit of a point $K$. Then let $\cB_K$ denote the \emph{greatest ball}
of centre $K$ inside the D-V cell $\cD(K)$, moreover let $\rho(K)$ denote the
\emph{radius} of $\cB_K$. It is easy to see that
$$
\rho(K)=\min_{\bg\in\Gamma\setminus\bI}\frac12 d(K,K^\bg).
$$
\end{definition}
The $\Gamma$-images of $\cB_K$ form a ball packing $\cB^\Gamma_K$ with centre
points $K^\bG$.
\begin{definition}\label{def:dens}
The \emph{density} of ball packing $\cB^\Gamma_K$ is
$$
\delta(K)=\frac{Vol(\cB_K)}{Vol\cD(K)}.
$$
\end{definition}
It is clear that the orbit $K^\Gamma$ and the ball packing $\cB^\Gamma_K$ have the
same symmetry group, moreover this group contains the starting
crystallographic group $\Gamma$:
$$
Sym K^\Gamma=Sym\cB^\Gamma_K\geq\Gamma.
$$
\begin{definition}\rm
We say that the orbit $K^\Gamma$ and the ball packing $\cB^\Gamma_K$ is
\emph{characteristic} if $Sym K^\Gamma=\Gamma$, else the orbit is not
characteristic.
\end{definition}
\subsection{Simply transitive ball packings}
\emph{Our problem is} to find a
point $K\in\ \SXR$ and the orbit $K^\Gamma$ for $\Gamma$ such that $\Gamma_K=\bI$
and the density $\delta(K)$ of the corresponding ball packing
$\cB^\Gamma(K)$ is maximal. In this case the ball packing $\cB^\Gamma(K)$ is
said to be \emph{optimal.}

The lattice of $\Gamma$ has a free parameter $p(\Gamma)$. Then we have to find the densest ball packing on $K$ for fixed
$p(\Gamma)$, and vary $p$ to get the optimal ball packing.
\begin{equation}
\delta^{opt}(\Gamma)=\max_{K, \ p(\Gamma)}(\delta(K)) \tag{3.1}
\end{equation}
Let $\Gamma$ be a fixed by {\it screw motions generated} group.
The stabiliser of $K$ is trivial i.e. we are looking the optimal kernel point
in a 3-dimensional region, inside of a fundamental domain of $\Gamma$ with free fibre parameter $p(\Gamma)$.
It can be assumed by the homogeneity of $\SXR$, that the fibre coordinate of the center of the optimal ball is zero.
\subsection{Optimal ball packings to space groups \textbf{1q.~I.}}
We consider those by screw motions generated groups in $\SXR$
whose point groups $\Gamma_0$ determine  spherical groups, characterized by their Macbeath signature.
The first set among them is {\bf 1q.~I.} which contain two space groups {\bf 1q.~I.~1}, {\bf 1q.~I.~2}:
\begin{gather*}
(+,0,[q,q],\{\})\times {1_\mathbf{R}}, ~ q\geq1, ~q\in \mathbf{N} \\
\Gamma_0=(g_1|g_1^q=1),
\end{gather*}
where $q$ is an integer parameter, which shows the degree of the rotation generating
the point group. It is easy to see that to the parameters $q=1$, $q=2$ do not exist ball packings thus we assume, that $3\le q \in \mathbf{N}$.
The possible translation parts of the generators of $\Gamma_0$ will be determined by the defining relations of the point group.
Finally, from the so-called Frobenius congruence relations we obtain two classes of non equivariant solutions:
$$\tau ~ \cong 0,~{\rm or} ~ ~ \frac{k}{q}, ~ k:=1 \dots \Big\lfloor \frac{q}{2} \Big\rfloor~ {\rm(i.e.~lower~ integer~ part~ of}~ \frac{q}{2}).$$
If $\tau \cong ~0$ then we get the $\SXR$ space group {\bf 1q.~I.~1} and if $\tau \cong ~
\frac{k}{q}, ~ k:=1 \dots \Big\lfloor \frac{q}{2} \Big\rfloor$ then we get the $\SXR$ space group {\bf 1q.~I.~2}.

We note here, that iff the greatest common divisor $(k,q)=1$, then the group $\Gamma= \mathbf{1q.~I.~2}$ will be fixed point free and the factor space
$\SXR/\mathbf{\Gamma}$ will be a compact orientable manifold (space form).

The fundamental domains of their point groups are spherical digons $B_0B_1$ with
angles $\frac{2\pi}{q}$ lying in the base plane $\Pi$ (see Fig.~1).

We shall apply the Cartesian homogeneous coordinate system introduced in Section 2.
\subsubsection{Optimal ball packing to space group \textbf{1q.~I.~1}}
We can assume, that $g_1$ is a rotation about the $z$-axis through angle $\frac{2\pi}{q}$ ($~q \ge 3$). The translation part in the
corresponding $\SXR$ space group is zero. Therefore, the group contain an arbitrary fibre translation $\tau$.
The fundamental domain ($D-V$ cell) of such a space group is a $\SXR$ prism with height $h=|\tau|$ whose volume can be calculated by the Theorem 2.3.

It can be assumed by the homogeneity of $\SXR$, that the fibre coordinate of the center of the optimal ball is zero and it is clear, that
the centre of the optimal ball is coincide with the centre of the digon $B_0C_0B_1C_1$ (see Fig.~1). Therefore,
the radius of the incircle or insphere of the optimal ball is $R^{opt}=\frac{\pi}{q}$. Moreover, we obtain the locally densest ball arrangement
if the height of the prism is $h_q^{opt}=2R_q^{opt}=\frac{2\pi}{q}$. These, locally densest ball arrangements, depending on $q$ are
denoted by $\mathcal{B}_q^{opt}$.
\begin{Theorem}
 The ball arrangement $\mathcal{B}_3^{opt}$ provides the optimal ball packing to the $\SXR$ space group $\mathbf{1q.~I.~1}$ with
 density $\delta^{opt}(\mathbf{1q.~I.~1})~\approx~0.50946$.
\end{Theorem}
{\bf{Proof}}

The fundamental domain of their point groups are a spherical digon $\mathcal{A}_q=B_0B_1$ with
angles $\frac{2\pi}{q}$ lying in the base plane $\Pi$ (see Fig.~1). The area of the a digon is $Vol(\cA_{q})=\frac{4\pi}{q}$,
while the radius of the incircle or insphere is $R^{opt}=\frac{\pi}{q}$. The volume of the locally optimal $D-V$ cell $\mathcal{P}_q$
(prism-like domain)
of the ball arrangement $\cB_q$ is $Vol(\mathcal{P}_q)=\frac{4\pi^2}{q^2}$ (see Theorem 2.3).
The volume $Vol(B_q^{opt})$ of the corresponding optimal ball $B_q^{opt}$ of radius $R^{opt}=\frac{\pi}{q}$ can be computed by the Theorem 2.2
and thus we obtain the densities $\delta^{opt}_q$ of the ball packings $\mathcal{B}_q^{opt}$
for arbitrary given parameter $3 \le q \in \mathbf{N}$. Moreover, if the parameter $q$ is sufficiently large then the volume of
the optimal ball $Vol(B_q^{opt})$ is approximately equal to the Euclidean one which is direct proportion to $\frac{1}{q^3}$.
Therefore, $\delta^{opt}(\mathbf{1q.~I.~1})$ is a strictly decrease function of $q$ and thus, $\mathcal{B}_3^{opt}$ provides
the densest ball arrangement. ~ ~$\square$
\begin{figure}[ht]
\centering
\includegraphics[width=7.5cm]{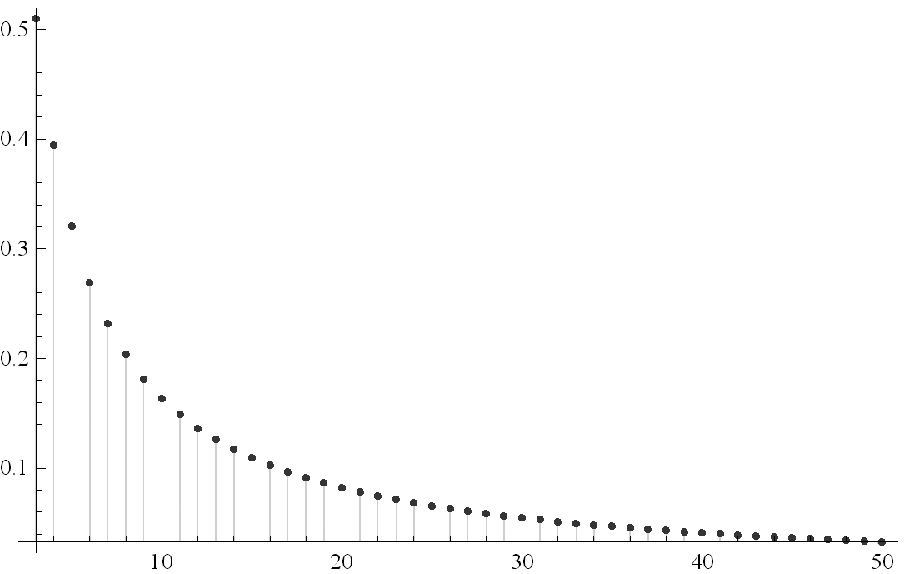} \includegraphics[width=6cm]{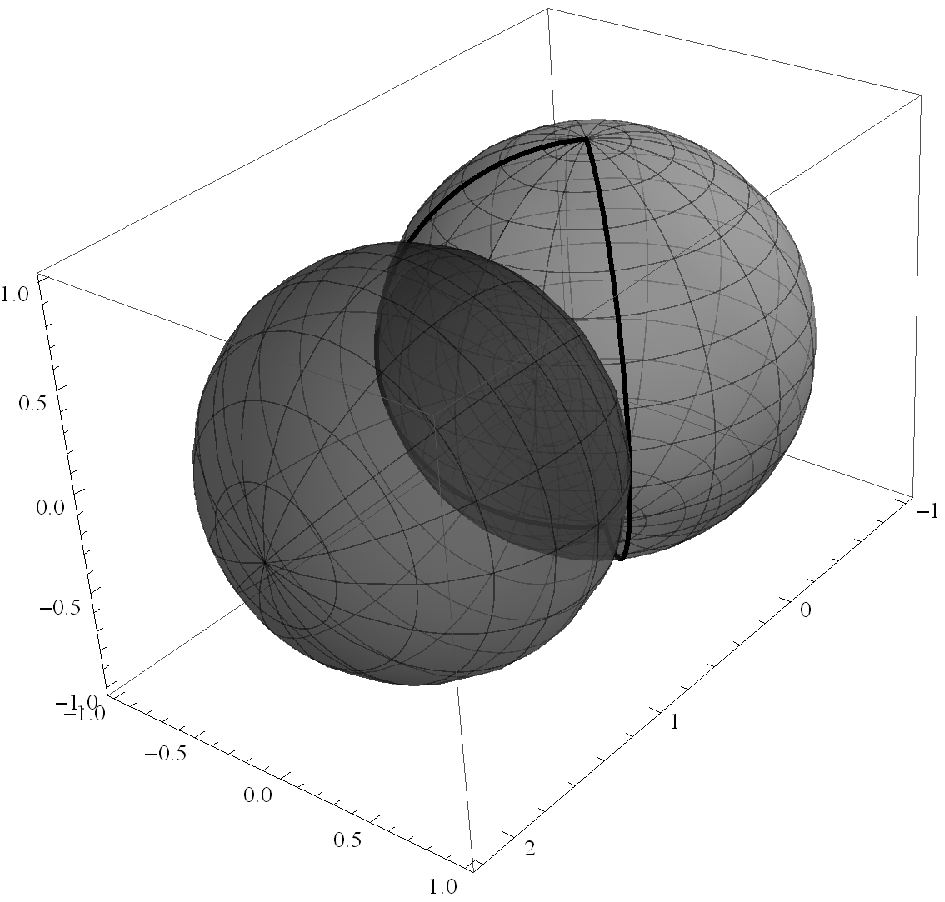}
\caption{The optimal densities of group $\mathbf{1q.~I.~1}$ for parameters $q=3,4,5,6,\dots$} and the optimal ball with the base plane $\Pi$ if $q=4$
\label{pic:1qdens}
\end{figure}

\subsection{Optimal ball packing to space groups \textbf{8.~I.}}\label{sec_8null}

We consider those screw motions generated groups \textbf{8.~I.~1} and  \textbf{8.~I.~2}
(using the Macbeath signature $(+,0,[2,3,3],\{\})\times {1_\mathbf{R}}$) whose point group $\Gamma_0$
is determined by the spherical group
\begin{gather*}
\Gamma_0=(g_1,g_2|g_1^2=g_2^3=(g_1g_2)^3=1).
\end{gather*}
where $g_1,g_2$ are two rotations with degrees of $2$ and $3$ respectively.
Fig.~3 shows a spherical fundamental domain $A_1A_2A_3A_4A_5$ of this group that is a spherical pentagon lying in the base plane
$\Pi$.
It is a spherical Dirichlet-Voronoi domain $\mathcal{D}(K)$ with kernel point $K$.
Moreover let $F$ be the midpoint of the spherical segment $A_1A_2$.
Here $g_1$ is the 2-rotation about the fibre line $f_F$ through the point $F$, $g_2$ is the 3-rotation about the fibre line $f_{A_3}$
through the point $A_3$ and $g_1g_2$ is the 3-rotation about the fibre line $f_{A_5}$, as well.

The possible translation parts of the generators will be determined by the defining relations of the point group.
From the Frobenius congruence relations we can obtain two classes of non-equivariant solutions:
$(\tau_1,\tau_2)\cong(0,0),\text{ or }(0,\frac{1}{3})$ and so two space groups \textbf{8.~I.~1} and \textbf{8.~I.~2}
(see \cite{F01}).
\subsubsection{Optimal ball packing to space group \textbf{8.~I.~1}}\label{sec:8i0}
First let us consider the space group $\Gamma_1=\mathbf{8.~I.~1}$, where the translation parts of the generators are $(\tau_1,\tau_2)=(0,0)$.
In this case the group is only generated by two rotations $g_1,g_2$ and the group contains a nonzero radial translation $\tau$.
The fundamental domain of the space group $\Gamma_1$ is a pentagonal prism $\mathcal{P}(K)=B_1B_2B_3B_4
B_5C_1C_2C_3C_4C_5$ which is derived from the spherical fundamental domain $\cD(K)=A_1A_2A_3A_4A_5$
by translations $\tau/2$ and $- \tau /2$ (see Fig.~3).
$\mathcal{P}(K)$ is also a $D-V$ cell of the considered group with kernel point $K$, as well.
\begin{figure}[h]
\centering
\includegraphics[height=6cm]{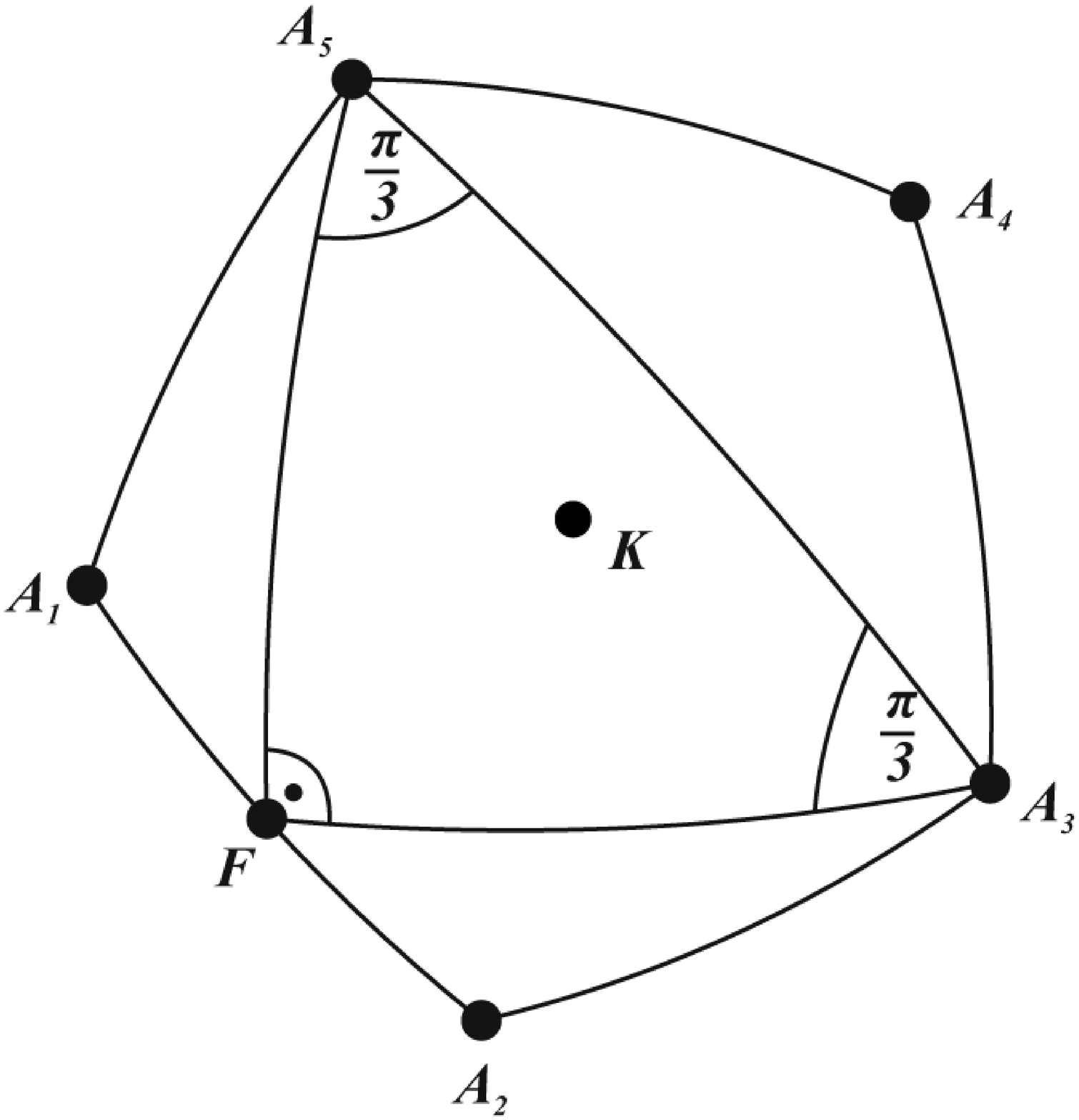} \includegraphics[height=6.5cm]{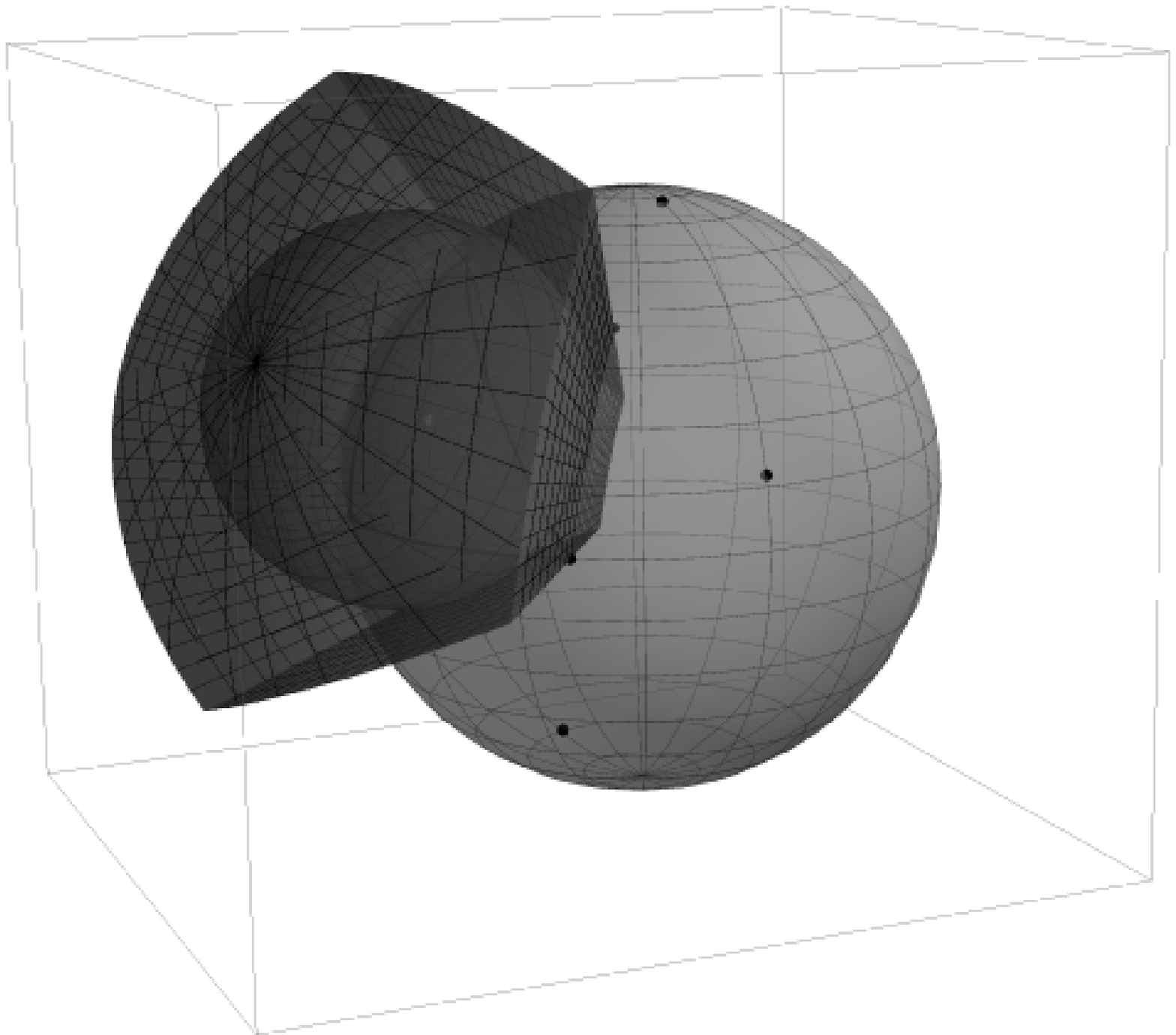}
\caption{The fundamental domain of the point group of the space group $\mathbf{8.~I.~1}$ and the optimal 5-gonal prism-like $D-V$-cell $\cP(K^{opt})$
with the base plane}
\label{pic:8i1}
\end{figure}
Let $\mathcal{B}(R)$ denote a geodesic ball packing of $\SXR$ space with balls $B(R)$ of radius $R$ where their
centres give rise to the orbit $K^{\Gamma_1}$. In the following we consider to each ball packing the {\it possible smallest
translation part} $\tau(K,R)$ (see Fig.~3) depending on $\Gamma_1$, $K$ and $R$.
A fundamental domain of $\Gamma$ is its prism-like $D-V$ cell $\cP(K)$ around the kernel point $K$.
It is clear that the optimal ball $\mathcal{B}_K$ has to touch some faces of its $D-V$ cell.
The volume of $\cP(K)$ can be calculated by the area of the spherical fundamental domain $\cD(K)$ and by the height $\tau(R,K)$.
The images of $\cD(K)$ form a congruent prism tiling by the discrete isometry group $\Gamma$.
For the density of the packing it is sufficient to relate the volume of the optimal ball
to that of the solid $\cP(K)$ (see Definition 3.4).
It is clear, that it is sufficient to study the kernel point as the point of the spherical triangle (lying in the base plane $\Pi$)
$FA_3A_5$ with angles $\frac{\pi}{3}$, $\frac{\pi}{2}$, $\frac{\pi}{3}$ without its vertices.

We shall apply the Cartesian homogeneous coordinate system introduced in Section 2.
\begin{Theorem}
The coordinates of the kernel point $K^{opt}$, the radius and the density of the optimal ball packings $\mathcal{B}^{opt}(R)$ of the
$\SXR$ space group $\mathbf{8.~I.~1}$ are
$$
K^{opt}=\Big(\sqrt{\frac{5+\sqrt{5}}{10}},\frac{1}{2}\sqrt{\frac{\sqrt{5}-1}
{\sqrt{5}}}, \frac{1}{2}\sqrt{\frac{\sqrt{5}-1}{\sqrt{5}}}\Big),
\ \ R^{opt}=\arccos\sqrt{\frac{5+\sqrt{5}}{10}}.
$$
$$
\delta^{opt}(\mathbf{8.~I.~1})=\frac{Vol(\mathcal{B}(R^{opt}))}{Vol(\cP(K^{opt}))}\approx 0.6005
$$
\end{Theorem}
{\bf{Proof}}

\textit{Our goal is to find the optimal kernel point $K^{opt}$ lying in the spherical triangle $FA_3A_5$, $K \notin \{F, A_3, A_5\}$
of the optimal Dirichlet-Voronoi cell $\cP(K^{opt})$ and the optimal
sphere radius $R^{opt}$, so that the density $\delta(K)$ of the ball packing is maximal.}
We can assume that $g_1$ is the 2-rotation about the $x$ axis, while the axis of 3-rotation $g_2$ is in the $[x,y]$ plane.

Because the translation parts of the generators $g_i$ $(i=1,2)$ in the corresponding $\SXR$ space group are zero,
therefore in this case the balls have a spherical shell-like arrangement $\cD(K)=A_1A_2A_3A_4A_5$
and the group contains a fibre translation $\tau$.

The optimal ball arrangement $\cB^{opt}(R^{opt})$ has to satisfy the following requirements:
\begin{enumerate}
	\item $\frac{\sqrt{3}}{2}d(K,K^{g_2})=d(K,K^{g_1})$
	\item $2R^{opt}=d(K,K^{g_1})=d(K,K^{\tau})$
\end{enumerate}
The first condition ensures the touching of the optimal sphere with the side neighbouring spheres of the prism $\cP$ on the level of the base
sphere, while the third condition requires the touching of the optimal sphere and its radial translated images.
\begin{figure}[h]
\centering
\includegraphics[height=11cm]{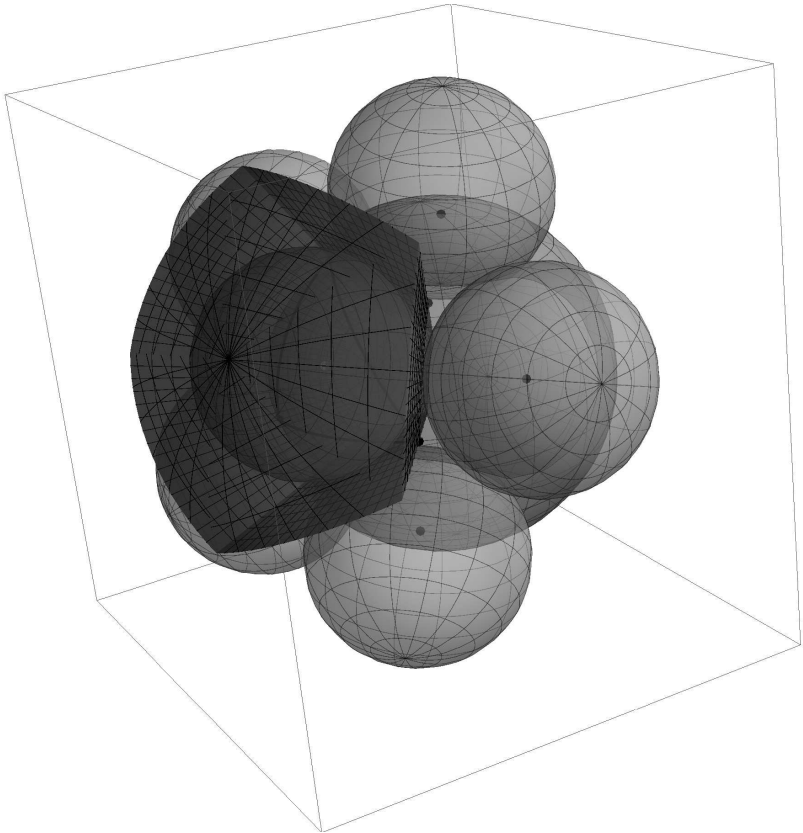}
\caption{The optimal ball arrangement $\cB^{opt}(R^{opt})$ with a optimal prism-like $D-V$ cell $\cP(K^{opt})$ of space group $\mathbf{8.~I.~1}$.}
\label{pic:8i1}
\end{figure}
By solving the above system of equations, we get the center $K^{opt}$ of the insphere
$K^{opt}=\Bigg(\sqrt{\frac{5+\sqrt{5}}{10}},\frac{1}{2}\sqrt{\frac{\sqrt{5}-1}
{\sqrt{5}}}, \frac{1}{2}\sqrt{\frac{\sqrt{5}-1}{\sqrt{5}}}\Bigg)$, and its radius $R^{opt}=\arccos\sqrt{\frac{5+\sqrt{5}}{10}}$.

It is easy to see, that the area of the base polygon $Vol(\cP(K^{opt}))=2\cdot Vol(FA_3A_5 \bigtriangleup)=\frac{\pi}{3}$,
so the volume of the Dirichlet-Voronoi cell can be computed by the Theorem 2.3.
Moreover, we get by Theorem 2.2 the volume of the insphere
$Vol(B(R^{opt}))\approx 0.6962$ and thus using the density formula given in Definition 3.4 we obtain the optimal density
$\delta^{opt}(\mathbf{8.~I.~1})\approx 0.6005$ (see Fig.~4). $\square$
\subsubsection{Optimal ball packing to space group \textbf{8.~I.~2}}
We get the space group $\Gamma_2 = \mathbf{8.~I.~2}$ if the translation part of the generators are
$(\tau_1,\tau_2)=(0,\frac{1}{3})$ i.e. $\tau_1=\tau;~\tau_2=\frac{\tau}{3}$.

Similarly to the above case the spherical fundamental domain of its point group is a spherical pentagon $\cD(K)=A_1A_2A_3A_4A_5$
lying in the base plane $\Pi$ (see Fig.~2.~a).
It can be assumed by the homogeneity of $\SXR$, that the fibre coordinate of the center of the optimal ball is zero
and it lies in the spherical triangle $FA_3A_5$, $K \notin \{F, A_3, A_5\}$ using the above denotations.

We shall apply for the computations the Cartesian homogeneous coordinate system introduced in Section 2.

We consider an arbitrary point $K(x^0,x^1,x^2,x^3)$ of spherical triangle $FA_3A_5$, $K \notin \{F, A_3, A_5\}$
in the above coordinate system in our model.

Let $\mathcal{B}(R)$ denote a geodesic ball packing of $\SXR$ space with balls $B(R)$ of radius $R$ where their
centres give rise to the orbit $K^{\Gamma_2}$. In the following we consider to each ball packing {\it the possible smallest
translation part} $\tau(K,R)$ depending on $\Gamma_2$, $K$ and $R$.
A fundamental domain of $\Gamma_2$ is its $D-V$ cell $\cP(K)$ with kernel point $K$ which is not a prism.
The optimal ball has to touch
some faces of its $D-V$ cell so that {\it the balls of the packing form a locally stable arrangement}.

The volume of $\cP(K)$ is equal to the volume of the prism which is given by the fundamental domain of the point group
$\Gamma_0$ of $\Gamma$ and by the height $|\tau(R,K)|$. The images of $\cP(K)$
by our discrete isometry group $\Gamma_2$ covers the $\SXR$ space without overlap.
For the density of the packing it is sufficient to relate the volume of the optimal ball
to that of the solid $\cP(K)$ (see Definition 3.4).
\begin{figure}[h]
\centering
\includegraphics[height=12cm]{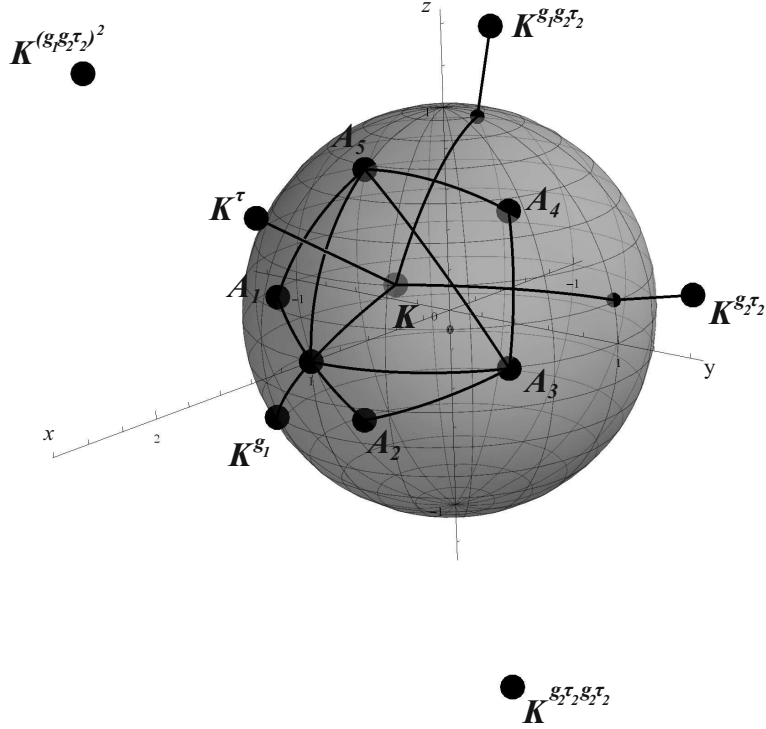}
\caption{A part of the orbit $K^{\Gamma_2}$ where $K(x^0,x^1,x^2,x^3) \in FA_3A_5$, $K \notin \{F, A_3, A_5\}$.}
\label{pic:8i1}
\end{figure}
It is clear, that the densest ball arrangement $\cB^{opt}_{\Gamma_2}(R)$ of balls $B(R)$ has to hold the following requirements:

\begin{enumerate}\label{req}
	\item $d(K,K^{g_1}) = 2 R = d(K,K^{g_2\tau_2})=d(K,K^{g_1g_2\tau_2})$
	\item $d(K,K^\tau) \geq 2 R$
	\item $d(K^{g_2g_2\tau_2\tau_2},K^{g_1g_2g_1g_2\tau_2\tau_2}) \geq 2 R$,
\end{enumerate}
where $d$ is the distance function of $\SXR$ space.

The ball packings described by the requirements above form a one parameter class of ball arrangements by the parameter $\tau$.
We can now examine the density function $\delta(\tau)$ and find its maximum.

We consider two main ball arrangements:

\begin{itemize}
	\item[a.] We denote by $\cB^{a}_{\Gamma_2}(R_a,K_a)$ the packing, where the above requirements and $d(K,K^\tau) = 2 R$ hold.
	(This is the case, where the kissing number $\Omega$ of the balls is maximal $\Omega=7$.)
	\item[b.] We denote by $\cB^{b}_{\Gamma_2}(R_b,K_b)$ the packing, where the above requirements and
	$d(K^{g_2g_2\tau_2\tau_2},K^{g_1g_2g_1g_2\tau_2\tau_2}) = 2 R$ hold.
	(This is the case, where the radius $R$ of the spheres is maximal, here $\Omega=5$)
\end{itemize}
\begin{Theorem}
 The ball packing $\cB^{a}_{\Gamma_2}(R_a,K_a)$ provides the optimal ball packing to the $\SXR$ space group $\mathbf{8.~I.~2}$ with
 density $\delta^{opt}(\mathbf{8.~I.~2})~\approx~0.6587$.
\end{Theorem}
{\bf{Proof}}

The corresponding $\tau$ values to the ball arrangements $\cB^{a}_{\Gamma_2}(R_a,K_a)$ and \linebreak $\cB^{b}_{\Gamma_2}(R_b,K_b)$
are $\tau_a \approx1.1608$ and $\tau_b \approx 3.5071$.
It is easy to see, that for any $\tau<\tau_a$ then $\delta(\tau) < \delta(\tau_a)\approx 0.6587$, and we get also that if $\tau=\tau_b$ then
the radius of the sphere is maximal, so for $\tau > \tau_b$ the density must be smaller than $\delta(\tau_b) \approx 0.5289$.

The derivatives of the density function $\delta(\tau)=\frac{Vol(B(\tau))}{Vol(\cP(K^{\tau}))}$ can be approximated by numerical methods to an arbitrary precision. Careful approximation and investigation of the first and second derivatives of the density function shows, that the second derivative is strictly positive on the interval $[\tau_a,\tau_b]$, therefore the maximum of the density function must be in one of the endpoints of the interval.

First we determine the coordinates of the points $K_i$, the radii $R_i$ $(i=a,~b)$ of the balls, the volumes of the corresponding balls and the
densities in both cases. Finally we get the following solutions, where the computations
were carried out by \textit{Wolfram Mathematica 8}:
\begin{gather*}
 R_a=0.5804, ~ K_a=(0.8362, 0.3877, 0.3877), ~  \delta(\tau_a) \approx 0.6587,\\
 R_b=0.7847, ~ K_b=(0.7075, 0.4996, 0.4996), ~  \delta(\tau_b) \approx 0.5289,
\end{gather*}

Therefore, the ball packing $\cB^{a}_{\Gamma_2}(R_a,K_a)$ provides the maximal density to the $\SXR$ space group $\mathbf{8.~I.~2}$. ~$\square$
\begin{figure}[h]
\centering
\includegraphics[width=10cm]{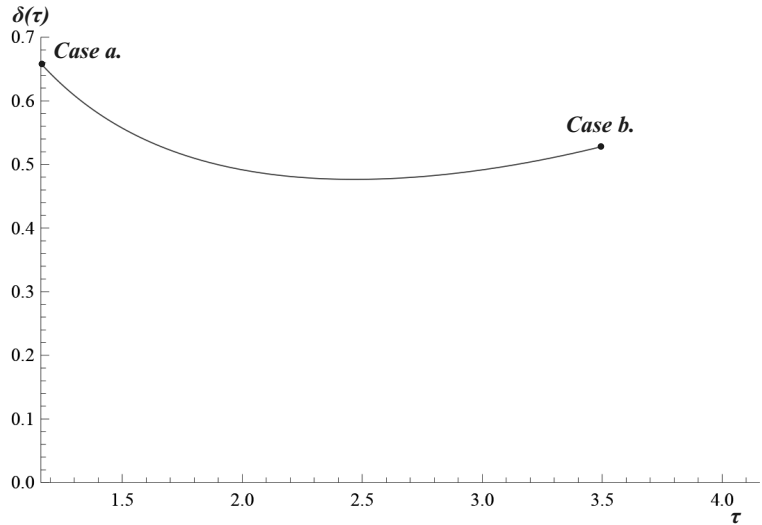}
\caption{The density function $\delta(\tau)$ of group \textbf{8.~I.~2}.}
\label{pic:dens}
\end{figure}

The remaining groups $\mathbf{1q.~I.~2},$ $\mathbf{3q.~I.~1}$, $\mathbf{3q.~I.~2}$, $\mathbf{3qe.~I.~3}$, $\mathbf{9.~I.~1}$, $\mathbf{9.~I.~2}$, $\mathbf{10.~I.~1}$
can be studied with the same way, as we have described previously. Finally we have summarized our results in the following Table.
Finally we get the following
\begin{Theorem}
The densest geodesic ball packings to $\mathbf{S}^2\!\times\!\mathbf{R}$ space groups generated by screw motions is derived
by the space group \textbf{3qe.~I.~3} with maximal packing density $\delta^{opt}(\mathbf{3qe.~I.~3}) \approx 0.7278$.
\end{Theorem}

\begin{table}[ht!]\centering
\begin{tabular}{|c|c|c|c|}
\hline
\textbf{Group}  & $(\tau_1,\tau_2)$ & Optimal radius & \textbf{Optimal density}\\
\hline
\hline
\textbf{1q.~I.~1} & $(0)$ & $R=\frac{\pi}{3}$  & $\delta\approx0.5094$ (with $q=3$)\\
\textbf{1q.~I.~2} & $(\frac{k}{q})$ & $R\approx1.1107$  & $\delta\approx0.5678$ (with $q=3,k=1$)\\
\textbf{3q.~I.~1} & $(0,0)$ & $R=\frac{\pi}{4}$ & $\delta\approx0.5919$ (with $q=3$)\\
\textbf{3q.~I.~2} & $(\frac{1}{2},\frac{1}{2})$ & $R\approx0.8417$ & $\delta\approx0.6758$\\
\textbf{3qe.~I.~3} & $(0,\frac{1}{2})$ & $R\approx0.8752$ & $\delta\approx0.7278$\\
\textbf{8.~I.~1} & $(0,0)$ & $R=\arccos\sqrt{\frac{5+\sqrt{5}}{10}}$ & $\delta\approx0.6004$\\
\textbf{8.~I.~2} & $(0,\frac{1}{3})$ & $R\approx0.5804$ & $\delta\approx0.6587$\\
\textbf{9.~I.~1} & $(0,0)$ & $R\approx 0.3812$ & $\delta\approx0.5758$\\
\textbf{9.~I.~2} & $(\frac{1}{2},0)$ & $R\approx0.4189$ & $\delta\approx0.6937$\\
\textbf{10.~I.~1} & $(0,0)$ & $R\approx0.2341$ & $\delta\approx0.5458$\\
\hline
\end{tabular}\label{table}
\end{table}

\end{document}